\documentclass[12pt]{amsart}

\usepackage[a4paper, hmargin=2.5cm, vmargin={3cm, 3cm}]{geometry}

\usepackage{amssymb, amsthm, amsmath}
\usepackage{xcolor}
\usepackage{graphicx}

\newcommand*\mcupinn[2]{\vcenter{\hbox{$\mathsurround=0pt
  \ifx\displaystyle#1\textstyle\else#1\fi\bigcup$}}}
\newtheorem{theorem}{Theorem}[section]

\DeclareGraphicsExtensions{.pdf,.png,.jpg}
\usepackage{pict2e}

\DeclareFontFamily{U} {MnSymbolC}{}

\DeclareFontShape{U}{MnSymbolC}{m}{n}{
  <-6> MnSymbolC5
  <6-7> MnSymbolC6
  <7-8> MnSymbolC7
  <8-9> MnSymbolC8
  <9-10> MnSymbolC9
  <10-12> MnSymbolC10
  <12-> MnSymbolC12}{}
\DeclareFontShape{U}{MnSymbolC}{b}{n}{
  <-6> MnSymbolC-Bold5
  <6-7> MnSymbolC-Bold6
  <7-8> MnSymbolC-Bold7
  <8-9> MnSymbolC-Bold8
  <9-10> MnSymbolC-Bold9
  <10-12> MnSymbolC-Bold10
  <12-> MnSymbolC-Bold12}{}

\DeclareSymbolFont{MnSyC} {U} {MnSymbolC}{m}{n}
\DeclareMathSymbol{\diamondtimes}{\mathbin}{MnSyC}{125}

\newtheorem*{theorem*}{Theorem}{\bf}{\it}
\newtheorem*{proposition*}{Proposition}{\bf}{\it}
\newtheorem*{observation*}{Observation}{\bf}{\it}
\newtheorem*{lemma*}{Lemma}{\bf}{\it}

\theoremstyle{definition}

\theoremstyle{remark}

\newcommand{\F}{\mathcal F}

\newcommand{\R}{\mathbb R}

\newcommand{\B}{\mathbb B}

\newcommand{\T}{\mathbb T}
\newcommand{\Sp}{\mathbb S}
\newcommand{\dv}{{\rm{div}}}
\newcommand{\dist}{{\rm{dist}}}
\newcommand{\diam}{{\rm{diam}}}

\newcommand{\h}{\mathcal{H}}

\def\XXint#1#2#3{{\setbox0=\hbox{$#1{#2#3}{\int}$ }
\vcenter{\hbox{$#2#3$ }}\kern-.6\wd0}}

\begin{document}
\title[]{Review of Yau's conjecture on zero sets of Laplace eigenfunctions}

\author[A. Logunov]{Alexander Logunov}
\address{A.L.: Department of Mathematics, Princeton University, Princeton, NJ, 08544;}
\email{log239@yandex.ru}
\thanks{This work was
completed during the time A.L. served as a Clay Research Fellow.}

\author[E. Malinnikova]{Eugenia Malinnikova}
\address{E.M.: Institute for Advanced Study, 1 Einstein Dr, Princeton, NJ 08540}

  \address{ Department of Mathematical Sciences,
Norwegian University of Science and Technology
7491, Trondheim, Norway}
\email{eugenia.malinnikova@ntnu.no}
\thanks{E.M. is supported  by  Project 275113
of the Research Council of Norway and NSF grant no. DMS-1638352}


\begin{abstract}  
This is a review of old and new results and methods related to the Yau conjecture on the zero set of Laplace eigenfunctions.
 The review accompanies two lectures given at the  conference CDM 2018.  
We discuss the works of Donnelly and Fefferman including their solution of the conjecture in the case of real-analytic Riemannian manifolds. 
 The review exposes the new results for Yau's conjecture in the smooth setting. We try to avoid technical details and emphasize the main ideas of the proof of Nadirashvili's conjecture. We also discuss two-dimensional methods to study zero sets.

 \end{abstract}
\maketitle

\nocite{*}

\setcounter{tocdepth}{1}
\tableofcontents
\section{Yau's conjecture}

  Yau conjectured \cite{Y} that for any $n$-dimensional $C^\infty$-smooth closed Riemannian manifold $M$ (compact and without boundary) 
 the Laplace eigenfunctions $\varphi_\lambda$ on $M$: $$\Delta \varphi_\lambda+ \lambda \varphi_\lambda =0,$$ satisfy

 $$c\sqrt \lambda \leq \h^{n-1}(\{\varphi_\lambda=0\}) \leq C\sqrt \lambda,$$  
  where $c$, $C$ depend only on the Riemannian metric on $M$ and are independent of the eigenvalue $\lambda$. The symbol $\h^{k}$ denotes the $k$ dimensional Hausdorff measure.


 The question of Yau is  connected to the quasi-symmetry conjecture, which states that 
 $$ c< \frac{\h^{n}(\{\varphi_\lambda>0\})}{\h^{n}(\{\varphi_\lambda<0\})}< C$$
for any non-constant eigenfunction $\varphi_\lambda$. In dimension two Yau's conjecture implies the quasi-symmetry conjecture, see the discussion  in Section \ref{sec: distr}.

The list of topics on geometry of Laplace eigenfunctions covered in this review is very limited. In particular we do not discuss the celebrated Courant nodal domain theorem, the variational methods, random eigenfunctions, the Kac-Rice formula and quantum ergodicity. The focus of this review is on the results and methods related to Yau's conjecture and the growth properties 
of solutions to elliptic PDE.  We would like to formulate some of  the previous  results on Yau's conjecture  and the  quasi-symmetry conjecture.
\begin{itemize}
\item[\textbullet] Brunning 1978 (\cite{B78}), Yau: Lower bound is true for $n=2$. 
\item[\textbullet]Donnelly \& Fefferman 1988 (\cite{DF}): Yau's conjecture and  the quasi-symmetry conjecture are true for real analytic metrics. In particular the conjectures are true for the spherical harmonics. 
\item[\textbullet]
{Nadirashvili 1988 (\cite{N88})}:  $n=2$, $\h^{1}(\{\varphi_\lambda=0\}) \leq C {\lambda} \log \lambda$.
\item[\textbullet]
{Donnelly \& Fefferman 1990 (\cite{DF1}), Dong 1992 (\cite{D92})}: $n=2$, $\h^{1}(\{\varphi_\lambda=0\}) \leq C \lambda^{3/4}$.
\item[\textbullet]
{Hardt \&  Simon 1989 (\cite{HS})}: $n\geq 2$,  $\h^{n-1}(\{\varphi_\lambda=0\}) \leq C \lambda^{C \sqrt \lambda}$.
\item[\textbullet] Nazarov \&  Polterovich \&  Sodin 2005 (\cite{NPS}): $n=2$, local bounds for asymmetry of sign. If $\varphi_\lambda(x)=0$, then 
for any $r>0$
 $$  \frac{\h^{2}(\{\varphi_\lambda>0\} \cap B_r(x))}{\h^{2}(\{\varphi_\lambda<0\} \cap B_r(x))}< C \log \lambda \log\log \lambda, \quad \lambda >10.$$
\item[\textbullet] {Colding \& Minicozzi 2011 (\cite{CMII}), Sogge \& Zelditch 2011 (\cite{SZ11}), 2012 (\cite{SZ12}), Steinerberger 2014 (\cite{St14})}:
$$ c\lambda^{\frac{3-n}{4}}\leq \h^{n-1}(\{\varphi_\lambda=0\}).$$
 
\end{itemize}

 In Section \ref{real-analytic df} we discuss the breakthrough work \cite{DF} of Donnelly  and Fefferman, which brought many ideas to nodal geometry.   The solution of the real-analytic case of Yau's conjecture uses the idea of the holomorphic extension for eigenfunctions.  The works of Donnelly and Fefferman explained how one can apply the methods of complex and harmonic analysis to nodal sets in the case when the Reimannian metric is real-analytic. Their works also gave us a point of view (for the general case of smooth metrics) that the geometry of nodal sets is related to the growth properties of the eigenfunctions.

In Section \ref{lift and frequency} we discuss two other useful ideas for the study of nodal sets: the harmonic extension for eigenfunctions and a powerful monotonicity property for harmonic functions.
  The latter monotonicity property can be formulated in two different  forms. The first form is the three balls inequality, which means that harmonic functions satisfy some sort of logarithmic convexity property. Agmon \cite{Ag} noticed that a logarithmic convexity property holds for  harmonic functions in the Euclidean space, Landis \cite{L63} proved a version of the three balls inequality for solutions of elliptic PDE with variable coefficients.
   The second form involves the notion of frequency for harmonic functions, which was introduced by Almgren \cite{Al}. The frequency function in a fixed ball is a characteristic of growth of the harmonic function in this ball.  Almgren did a remarkable discovery that the frequency function is monotone with respect to the radius when the center of the ball is fixed. Garofalo and Lin \cite{GL} proved a version of the monotonicity property of the frequency for harmonic functions on smooth Riemannian manifolds, which has many applications for nodal sets.

In Sections \ref{upper bound} and \ref{lower bound} we expose the recent results in the smooth case including the polynomial upper bound \cite{L1} and the lower bound \cite{L2} in Yau's conjecture:
$$c\sqrt \lambda \leq \h^{n-1}(\{\varphi_\lambda=0\})  \leq C \lambda^{C_n}.$$ 
 The recent results follow the path suggested by Nadirashvili \cite{N97}, who argued that there is no hope to understand nodal sets if we don't understand zero sets of harmonic functions. In order to attack Yau's conjecture Nadirashvili formulated  two conjectures on harmonic functions in the three dimensional Euclidean space. One of them was recently solved \cite{L2} and the proof of Nadirashvili's conjecture implied the lower bound in Yau's conjecture. The second conjecture was asked by several mathematicians and goes back to at least Lipman Bers. We will formulate it in Section \ref{Cauchy problem}. The conjecture concerns the Cauchy uniqueness problem and is still  open. However weaker results on unique continuation properties of elliptic PDE were used to prove polynomial upper bounds \cite{L1} in Yau's conjecture. 

 Let us also mention that in dimension $n=2$,
 one can improve  $3/4$ from Donnelly-Fefferman's bound by a tiny $\varepsilon$ (\cite{LM}):
$$\h^{1}(\{\varphi_\lambda=0\}) \leq C \lambda^{3/4-\varepsilon}.$$
 The conjectured upper bound is still a challenging problem even in dimension two where a lot of tools are available. We describe two-dimensional methods in Section \ref{sec:2d}.

\section{Introduction: eigenfunctions, zeros and growth}
 
\subsection{Eigenvalues of Laplace operator}
 We briefly recall the basic properties of Laplace eigenfunctions and refer to \cite{Ch2},\cite{Ch1} for the introduction to the subject.
Let $M$ be a closed manifold (compact and without boundary) with a given Remannian metric $g$. We denote by $\Delta$ the Laplace operator on $M$ defined by this metric. An eigenfunction of the Laplace operator on $M$ is a solution of the equation
\[\Delta\varphi+\lambda \varphi=0.\]
The operator $-\Delta$ is non-negative and has a discrete spectrum,
\[0=\lambda_0<\lambda_1\le \lambda_2\le ...\ .\] 
The smallest eigenvalue is $\lambda_0=0$ and the corresponding eigenfunction is constant. Eigenfunctions that correspond to distinct eigenvalues are orthogonal: $\int_M \varphi_k \varphi_l =0$.

 For a subdomain  $\Omega$ of $M$  with (piecewise)  smooth boundary the eigenfunctions of the  Laplace operator on $\Omega$ with Dirichlet boundary conditions are solutions of the problem
\[\begin{cases}
\Delta\varphi+\lambda\varphi=0\quad{\text{in}}\ \Omega\\
\varphi=0\quad{\text{on}}\ \partial\Omega.
\end{cases}  \]
All eigenvalues of the Dirichlet Laplacian are positive,
\[0<\lambda_1(\Omega)<\lambda_2(\Omega)\le...\ .\]
The first eigenvalue is simple and the corresponding first eigenfunction does not change sign in $\Omega$. The eigenfunctions corresponding to the higher eigenvalues are orthogonal to the first one and take both positive and negative values in $\Omega$. 
There is a variational characterization of the eigenvalues known as the Rayleigh quotient. The first eigenvalue is given by 
\[\lambda_1(\Omega)=\inf_{f}\frac{\int_{\Omega}|\nabla f|^2}{\int_{\Omega}|f|^2},\]
where the infimum is taken over all non-zero functions $f\in C^1(\bar{\Omega})$ such that $f=0$ on $\partial \Omega$. This implies in particular that if $\Omega_0\subset\Omega$ then
\[\lambda_1(\Omega_0)\ge\lambda_1(\Omega).\]

We denote by $j_n$ the first eigenvalue of the Dirichlet Laplace operator for the unit ball $\B^n\subset\R^n$. Then a simple renormalization implies that $\lambda_1(B_r)=r^{-2}j_n$ for the $n$-dimensional Euclidian ball of radius $r$. If $M$ is a closed Riemannian $n$-dimensional manifold, then (using the Rayleigh quotient and the normal coordinates) one can check that 
\[\lim_{r\to 0}r^2\lambda_1(B_r(x))=j_n\] 
for any $x\in M$, where $B_r(x)$ is the geodesic ball centered at $x$; the limit is uniform in $x$, see \cite[chapter 3.9]{Ch2}.

\subsection{Density of zeros.} \label{density}
Suppose that $\varphi_\lambda$ is a non-constant eigenfunction of the Laplace operator on $M$,  then $\int_M \varphi_\lambda=0$ and it changes sign. We consider the zero set of $\varphi_\lambda$,
\[Z(\varphi_\lambda)=\{x\in M: \varphi_\lambda(x)=0\},\]
and the connected components of its complement, $M\setminus Z(\varphi_\lambda)=\cup_j\Omega_j$. The domains $\Omega_j$ are called nodal domains of the eigenfunction $\varphi_\lambda$. The restriction of $\varphi_\lambda$ onto each domain $\Omega_j$ is an eigenfunction of the Dirichlet Laplace operator and, since $\varphi_\lambda$ does not change sign in $\Omega_j$, it is the first eigenfunction (we skip the discussion of the regularity properties of the domains,  one can find the details in \cite{Che} and \cite[chapter 1.5]{Ch1}). Therefore $\lambda_1(\Omega_j)=\lambda$ for each $\Omega_j$.
 
Now it is easy to see that $Z(\varphi_\lambda)$ is $c/\sqrt{\lambda}$ dense in $M$. If $x\in M$ and $\dist(x, Z(\varphi_\lambda))>r$ then $B_r(x)\subset\Omega_j$ for some $j$. It implies that
\[\lambda=\lambda_1(\Omega_j)<\lambda_1(B_r(x))<C(M)r^{-2}.\]
Hence $r<c/\sqrt{\lambda}$.

The lower bound in Yau's conjecture is supported by the fact that the zero set of $\varphi_\lambda$ is $\frac{C}{\sqrt \lambda}$ dense on $M$.

\subsection{Two  examples.} Let $\T^n$ denote the standard torus. We identify it with the cube $[0,1]^n$ with glued opposite faces. There is a basis for $L^2(\T^n)$ consisting of eigenfunctions of the Laplace operator. The elements are products of trigonometric functions (sines and cosines) with frequencies that are integer multiples of $2\pi$. For example
\[\phi(x_1,...,x_n)=\sin(2\pi k_1x_1)...\sin(2\pi k_nx_n)\]
is an eigenfunction with the eigenvalue $\lambda=4\pi^2(k_1^2+...+k_n^2)$. The zero sets of such eigenfunctions considered on the cube $[0,1]^n$ are unions of hyperplanes parallel to the coordinate hyperplanes.

Another  example is the unit sphere $\Sp^n$. The eigenfunctions on $\Sp^n$ are restrictions of homogeneous harmonic polynomials in $\R^{n+1}$. For $n=2$ we get the classical spherical harmonics. There is a basis consisting of spherical harmonics, whose zero sets are unions of "meridians" and circles of constant "latitude".

In both examples the zeros sets for basis eigenfunctions look very regular. However for the torus and for the sphere of dimension larger than one the multiplicities of the eigenvalues of the Laplace operator can be large. Interesting examples appear when we take zero sets of linear combinations of basic eigenfunctons corresponding to the same eigenvalue. 
A beautiful topic that we do not discuss here is the zero sets of random eigenfunctions (linear combinations with random coefficients). The interested reader can start wtih   \cite{NS2010},\cite{RW08},\cite{Can2019},\cite{Ze2009},\cite{KKW} for the introduction to random eigenfunctions.

\subsection{Vanishing order.}
 
 It is said that  the vanishing order of a smooth function $f$ at a point $x$ is $k$ if any derivative of order smaller than $k$
of $f$ at $x$ is zero and there is some non-zero derivative of order $k$. The vanishing order of $f$ at $x$ is zero if $f(x)\neq 0$
 and $\infty$ if all derivatives of any order at $x$ are zero.

In dimension two the vanishing order of any  eigenfunction  $\varphi_\lambda$ at a point $x$ has a geometrical meaning of the number of 
 nodal curves intersecting at $x$ and the nodal curves have equiangular intersection  at this point \cite{B55}. 
 A very natural question is
 how large could be the vanishing order of $\varphi_\lambda$?

 Donnelly and Fefferman \cite{DF} answered  this question for smooth Riemannian manifolds of any dimension  by showing that vanishing order of $\varphi_\lambda$ at any point is smaller
 than $C\sqrt \lambda$. The result is sharp if we don't make any extra assumptions on the Riemannian manifold. There are spherical harmonics with  vanishing order comparable to $\sqrt \lambda$.  Peter Sarnak suggested that for surfaces with negative curvature the bound on the vanishing order should be improved to $c_\varepsilon \lambda^\varepsilon$ for any $\varepsilon>0$, and it would imply a good bound on the multiplicity of the eigenvalues.
 
 The proof \cite{DF} of the doubling index  estimate $C\sqrt \lambda$ for general closed Riemannian manifolds is using Carleman inequalities and is inspired by the paper \cite{Ar} of Aronszajn on unique continuation properties of solutions of elliptic PDE of second order. The basic
question of unique continuation is whether a non-zero solution can vanish on an open set.

\subsection{Doubling index.}\label{ss:di}
 Donnelly and Fefferman also proved \cite{DF} a useful bound for the growth of $\varphi_\lambda$:
 \begin{equation} \label{eq:di}
 \log \frac{\sup_{2B} |\varphi_\lambda|}{\sup_{B} |\varphi_\lambda|} \leq C_M \sqrt \lambda 
\end{equation} 
for any geodesic ball $B$ on $M$ and the geodesic ball  $2B$ with the same center and twice bigger radius than $B$.

 The number  $\log_2 \frac{\sup_{2B} |f|}{\sup_{B} |f|}$ is called the doubling index of the function $f$ in the ball $B$
and is denoted by $N_f(B)$. Note that for $C^\infty$ smooth functions
$$\lim\limits_{r\to 0}  N_f(B_r(x)) = \textup{ vanishing order of $f$ at $x$}.$$   
The doubling index estimate \eqref{eq:di} implies the bound $C \sqrt \lambda$ for the vanishing order for eigenfunctions.

\subsection{BMO norm of $\log|\varphi_\lambda|$}\label{ss:BMO}
 One of the ideas of the works of Donnelly and Fefferman is that the eigenfunctions $\varphi_\lambda$ behave as  polynomials of degree $\sqrt{\lambda}$. In particular they prove \cite{DF2} Bernstein type inequalities for the norms of $\nabla \varphi_\lambda$ and conjecture that $\|\log|\varphi_\lambda|\|_{BMO}\le C\sqrt{\lambda}$, see \cite{ST93} for the definition of BMO space. If $P$ is a polynomial of one variable of degree $d$ then $P(z)=a(z-z_1)..(z-z_d)$ and it is clear that $\|\log |P(z)|\|_{BMO}\le Cd$. Similar estimate holds for polynomials of several variables.

Donnelly and Fefferman showed that $\|\log|\varphi_\lambda|\|_{BMO}\le C\lambda^{n(n+2)/4}$. This estimate was improved \cite{CM} by Chanillo and Muckenhoupt  and then \cite{Lu} by Lu, and \cite{HanLu} by Han and Lu.  The recent result \cite{LM2},\cite{LM1} on quantitative unique continuation for solutions of second order elliptic PDEs  implies the conjectured  bound  $\|\log\varphi_\lambda\|_{BMO}\le C\sqrt{\lambda}$. The conjectured bound appeared to be connected to a question of Landis, which will be  discussed in Section \ref{Landis spheres}.

We  can rewrite the estimate for the BMO-norm of $\log|\varphi_\lambda|$  as a propagation of smallness result. If for some constant $c>0$, each cube $Q$, and each $a>0$
\[\h^n\{x\in Q:|\varphi_\lambda(x)|\le e^{-a}\sup_Q|\varphi_\lambda|\}\le Ce^{-ca/\sqrt{\lambda}}|Q|\]
then $\|\log|\varphi_\lambda|\|_{BMO}\le C\sqrt{\lambda}$. The last inequality is equivalent to the following estimate
\[\sup_{Q}|\varphi_\lambda|\le C\left(C\frac{|Q|}{|E|}\right)^{C\sqrt{\lambda}}\sup_E|\varphi_\lambda| ,\]
for any subset $E\subset Q$ with positive measure. Note that this inequality resembles the classical inequality of Remez \cite{Rem,B} for polynomials.

Looking at the spherical harmonics $u(x,y,z)= \Re (x+iy)^n$ one can see  that  Laplace eigenfunctions can be $e^{-c\sqrt \lambda}$ small on  a fixed open subset of the manifold. Such localization cannot happen on the standard torus $\T^n$. Various strong results on the torus with the standard metric were obtained \cite{BR11a},\cite{BR11b},\cite{BR12} by Bourgain and Rudnick. In particular they proved a uniform $L^2$ restriction bounds on curves for $2$-dimensional torus.

For negatively curved Riemannian manifolds we believe that it is possible to prove better versions of the BMO estimate and better bounds for the doubling index.
We would like to mention an outstanding recent result  by Bourgain \& Dyatlov \cite{BD} and Dyatlov \& Jin \cite{DJ}.
\begin{theorem}[\cite{BD},\cite{DJ}]
 Under assumption that $(M,g)$ is a closed Riemannian surface with constant negative curvature the following inequality holds for all Laplace eigenfunctions $\varphi_\lambda$ on $M$.
 For any open subset $E$ of $M$ there exists $c=c(E,M,g) >0$ (independent of the eigenvalue $\lambda$)  such that 
$$ \int_E \varphi_\lambda^2 \geq c  \int_M \varphi_\lambda^2.$$
\end{theorem}

\section{Zero sets of eigenfunctions on surfaces} \label{sec:2d}
\subsection{Local structure of the zero set in dimension two} Let $M$ be a surface with a given Riemannian metric. Locally at small scales the zero set of any eigenfunction on $M$ looks like the zero set of a harmonic function on the plane. If $\varphi_\lambda(x)=0$ and the vanishing order of $\varphi$ at $x$ is $k$ (the number of  derivatives of $\varphi$, which are zero at $x$), then in a small geodesic ball centered at $x$ the zero set $Z(\varphi_\lambda)$ consists of $k$ smooth curves intersecting at the point $x$ and forming equal angles $\pi/k$ at the point of intersection. To the best of our knowledge it was first observed by Bers \cite{B55}. In higher dimensions it is also true that the first term in the Taylor series of $\varphi_\lambda$ is a harmonic function, however the local structure of the zero set in this case can be quite complicated (even in the Euclidean space $\mathbb{R}^3$) and not stable, see an example of a harmonic polynomial in \cite{LM2015}.

\subsection{Estimates of the length of the zero set from below.} The lower bound in the Yau conjecture on surfaces was proved \cite{B78} by Bruning and also by Yau. The result follows from the density estimate of the zero set and an observation on the diameter of connected components of the zero set.

\textbf{Observation.} There exists a constant $c=c(M)$ such that if $\varphi_\lambda$ is an eigenfunction and $\Omega_j$ is a connected component of $M\setminus Z(\varphi_\lambda)$ then $\diam(\Omega_j)\ge c(M)/\sqrt{\lambda}$. To prove the observation, suppose that $\Omega_j\subset B_r$ for some geodesic ball $B_r$ with radius $r$. Then by the monotonicity of the first eigenvalue,
\[\lambda=\lambda_1(\Omega_j)\ge\lambda_1(B_r)> c_1r^{-2}.\]
It implies that $\diam(\Omega_j)\ge c_2\lambda^{-1/2}.$

 The latter observation implies that if $\varphi_\lambda(x)=0$, then $$\h^1(\{\varphi_\lambda =0\}\cap B_{1/\sqrt \lambda}(x)) \geq c/\sqrt \lambda.$$
Combined with the fact that $\{\varphi_\lambda =0\}$ is $C/\sqrt \lambda$ dense on $M$, the observation implies the lower bound in Yau's conjecture in dimension two.

\subsection{Singular points} Donnelly and Fefferman obtained a number of interesting results \cite{DF1} on the zero sets for eigenfunctions on surfaces. They also considered the singular set of  eigenfunctions, defined as
\[S(\varphi_\lambda)=\{x\in M: \varphi_\lambda(x)=0, \nabla\varphi_\lambda(x)=0\}.\]
For eigenfunctions on surfaces the singular set is a discrete set of points and the number  of singular points of $\varphi_\lambda$ is bounded by $C\lambda$.  Donnelly and Fefferman proved a stronger statement:\\
If $B$ is a geodesic ball on $M$ of radius $c\lambda^{-1/4}$, and for each point $p\in B\cap S(\varphi_\lambda)$ let the vanishing order of $\varphi_\lambda$ at $p$ be $k(p)+1$ (vanishing order at any singular point is at least two), then
\begin{equation}\label{eq:sp}
\sum_{p} k(p)\le c\sqrt{\lambda}.
\end{equation}
The estimate \eqref{eq:sp} is sharp in several ways: one cannot 
enlarge the radius $c\lambda^{-1/4}$,
 there exist spherical harmonics with vanishing order at one point comparable to $\sqrt{\lambda}$ and  there are also spherical harmonics with the total number of singular points comparable to $\lambda$. 

One of the tools used by Donnelly and Fefferman is a simple and powerful two dimensional Carleman inequality. Let $D$ be a domain on the complex plane and $h$ be a smooth function on $D$. For any complex valued $f\in C^\infty_0(D)$ the following inequality holds
\[\int_D|\bar\partial f|^2e^h\ge \frac{1}{4}\int_D\Delta h|f|^2e^h,\]
and if $\Delta h \geq  1 $ in $D$, then
\[\int_D|\Delta f|^2e^{th}\ge c{t^2}\int_D|f|^2e^{th}.\]

A remarkable idea \cite{DF1} due to Donnelly and Fefferman explains how two dimensional Carleman inequalities help to unite the information on the behavior of the eigenfunctions  (such as growth, vanishing order, doubling index) near several points. This method is quite flexible in dimension two.
The original Carleman approach \cite{C39} concerns the two dimensional case, but there are higher dimensional generalizations of Carleman inequalities (\cite{Ar},\cite{Co},\cite{T95}). They allow to study the behavior of a solution to elliptic PDE near one point or near infinity, 
but higher dimensional Carleman inequalities are less flexible. The conditions on Carleman weights in higher dimensions are hard to apply in the situations where you have to work with the behavior of the eigenfunctions near several points or curves.

In Section \ref{Cauchy problem}  we will formulate an old open question, which shows that we don't understand  unique continuation properties for elliptic PDE well enough in higher dimensions.

\subsection{Estimate of the zero set from above by Donnelly and Fefferman} The upper bound in the Yau conjecture for eigenfunctions on surfaces with smooth Riemannian metric is still an open problem. Donnelly and Fefferman showed \cite{DF1} that
\begin{equation}\label{eq:3/4}
\h^1(Z(\varphi_\lambda) )\le C\lambda^{3/4}.
\end{equation}
Once again they worked on the scale $c\lambda^{-1/4}$ and used the  bound $C\sqrt \lambda$ for the doubling index in any ball. They showed that $\h^1(Z(\varphi_\lambda)\cap B)\le C\lambda^{1/4}$ for any geodesic ball $B$ of radius $c\lambda^{-1/4}$ by proving the following estimate for solutions of elliptic inequalities.

 \textbf{Estimate for the length of nodal set in dimension two} (\cite{DF1}). Let $Q$ be the unit square and $N>1$.
{ Suppose that a function $
	\varphi:2Q\to \R$ satisfies
	$$|\Delta \varphi| \leq N |\varphi|$$ 
	in $2Q$ and  for any subcube $q$ of $Q$ the doubling index of $
	\varphi$ in $q$ is smaller than $N$. Then 
	\[\h^1\left(x\in \frac{1}{100}Q: 
	\varphi(x)=0\right)\le CN.\] }  
The proof of the latter statement is not simple. It is based on a two-dimensional Carleman inequality with a carefully chosen weight adjusted to the function $\varphi$ and the idea of the Calderon--Zygmund decomposition. 

As the result Donnelly and Fefferman proved that on the scale $\frac{1}{\sqrt \lambda}$ the length of zero set can be estimated from above in terms of the doubling index. It is remarkable that there is also a lower bound. 

\textbf{Estimate of the length of nodal lines in terms of the doubling index} (\cite{R-F}, \cite{DF1},\cite{NPS}).
\begin{equation} \label{eq: length vs growth}
 c N_{\varphi_\lambda}(B_{\frac{1}{4\sqrt\lambda}}(x))-C \leq \sqrt \lambda \cdot \h^1(\{\varphi_\lambda =0\} \cap B_{\frac{1}{\sqrt\lambda}}(x))  \leq C N_{\varphi_\lambda}(B_{\frac{2}{\sqrt\lambda}}(x))+C.
\end{equation}

\textbf{Remark.}   The recent combinatorial argument  \cite{LM} shows that there exists $\varepsilon>0$ such that for any closed surface $M$ the eigenfunctions $\varphi_\lambda$ on $M$ satisfy
\[\h^1(Z(\varphi_\lambda))\le C\lambda^{3/4-\varepsilon}\]
for some $C=C(M)$. It demonstrates that there is a room for improvement for estimate \eqref{eq:3/4}. A challenging problem is to prove the upper bound conjectured by Yau even in dimension two.

\subsection{Yau's conjecture and  distribution of doubling indices} \label{sec: distr}

 Let $M$ be covered by $\sim \lambda $ geodesic discs $B_i$ of radius $C/\sqrt \lambda$ so that each
 point of $M$ is covered at most $C_1$ times and $\varphi_\lambda$ is zero at the center of each disc $B_i$.

\textbf{Conjecture}(Nazarov,Polterovich and Sodin). There is a numerical constant $C$ (independent of $\lambda$ and of the covering) such that
 $$ \frac{\sum N(B_i)}{\# B_i}\leq C.$$

 In view of \eqref{eq: length vs growth} the latter conjecture is equivalent to the Yau conjecture in dimension $2$.

\textbf{Weak form of the conjecture.}
 At least half of $B_i$ have a bounded doubling index.

 \textbf{Comment.} The weak  conjecture implies the quasi-symmetry conjecture: $$c <\frac{\h^2(\varphi_\lambda > 0)}{\h^2(\varphi_\lambda < 0)} < C. $$

\subsection{Approach of Dong} A different method to study  zeroes and singular sets of eigenfunctions on surfaces was suggested by Dong \cite{D92}. Let  $\varphi_\lambda$ be a Laplace eigenfunction on $2$-dimensional manifold $M$. Dong considered the function $q=|\nabla \varphi_\lambda|^2+\lambda\varphi_\lambda^2/2$ and obtained an estimate for $\Delta \ln q$ outside of the singular set of $\varphi_\lambda$. For a harmonic function $u$ on the Euclidean plane we know that $\nabla u$ can be identified with an analytic function and $\ln |\nabla u|$ is subharmonic. This simple fact has a remarkable power in complex and harmonic analysis. Dong proved that
\[\Delta \ln q\ge-\lambda+2\min(K,0)\]
on the set $\{ q\neq 0 \}$, where $K$ is the Gaussian curvature of the surface. Using the latter inequality Dong found different proofs for the estimate \eqref{eq:sp} of the sum of the vanishing orders at the singular points  and for the Donnelly--Fefferman bound \eqref{eq:3/4} of the length of the nodal set.

\subsection{Applications of quasiconformal mappings to eigenfunctions.}
 On the scale $c/\sqrt \lambda$ the Laplace eigenfunction $\varphi_\lambda$ behaves like a harmonic function. For Riemannian surfaces one can justify the latter claim in a rigourous way involving quasiconformal mappings. We will briefly discuss the reduction, and the interested reader can read the details in \cite{NPS} and learn more about the quasiconformal mappings and their applications to PDE in \cite{AIG}.
 There are several steps in the reduction.

 For Riemannian surfaces it is convenient to work in local conformal coordinates: the equation for the eigenfunctions simplifies to   $$\Delta \varphi + \lambda q \varphi=0,$$ where $\Delta$ is the standard Laplace operator on $\R^2$ and $q$ is a bounded function. If we consider a disc of radius $\epsilon/\sqrt \lambda$ and rescale it to the unit disc, the equation reduces to
$$  \Delta \varphi + V \varphi=0$$
 in $D=\{|z|<1\}$ with $\|V\|_\infty<C\epsilon^2$.

\textbf{Claim.} If $\| V\|_\infty$ is sufficiently small, then there is a positive solution $f$ to the  equation $$  \Delta f + V f=0$$
in $D$ such that $$1 - C \| V\|_\infty \leq f\leq 1.$$

 The ratio $u= \varphi / f$ satisfies in $D$ the equation in divergence form:
 $$ \dv(f^2  \nabla u) =0.$$

 The theory of quasiconformal mappings joins the game here. There is
 a $K$-quasiconformal homeomorphism $g$ from $D$ to $D$ with $g(0)=0$
 such that  $h = g \circ u $ is a harmonic function in $D$ and  $K>1$ satisfies
 $$ \frac{K-1}{K+1} \leq \sup_D \frac{1-f^2}{1+f^2}.$$
In particular $K$ tends to 1 as $f$ becomes close to 1 in $D$.
 We do not know the change of variables  $g$ explicitly, moreover it depends on the auxiliary function $f$,
 but  $g$ possesses good geometric properties with quantitative estimates that depend only on $K$, which is under control. 
For instance,  Mori's theorem states that $g$  is $1/K$- H\"older 
and 
$$\frac{1}{16} |z_1-z_2|^K\leq |g(z_1) - g(z_2)| \leq 16 |z_1-z_2|^{1/K}.$$

Nadirashvili \cite{N88} suggested to apply quasiconformal mappings to get the bounds for the length of zero sets of eigenfunctions.
Nazarov, Polterovich, Sodin \cite{NPS} applied quasiconformal mappings and Astala's area distortion theorem to  asymmetry of sign of the eigenfunctions.
 They showed that  if $\varphi_\lambda(x)=0$, then 
for any geodesic disc $B_r(x)$
 $$  \frac{\h^{2}(\{\varphi_\lambda>0\} \cap B_r(x))}{\h^{2}(\{\varphi_\lambda<0\} \cap B_r(x))}< C \log \lambda \log\log \lambda.$$
Above we assume that $\lambda > 10$ to make $\log\log \lambda$ well defined.

To prove the bound above it sufficient to consider the case $r\leq c /\sqrt \lambda$ only, the case of bigger scales follows in a straightforward way using the fact that the nodal set is $C/\sqrt\lambda $ dense. Nazarov, Polterovich and Sodin used the quasiconformal mappings and  the doubling index bound $C\sqrt \lambda$ 
to reduce the  question of quasi-symmetry for eigenfunctions in dimension 2 to a question about harmonic functions on the plane.

\textbf{Quasi-symmetry of sign of harmonic functions with controlled growth.}
 Let $u$ be a harmonic functions in $\R^n$ with $u(0)=0$. Suppose that $N_u(B) \leq N$.
 How large $\frac{|\{u>0 \}\cap B|}{|\{u<0 \}\cap B|}$ can be?

Nazarov, Polterovich, Sodin answered this question in dimension $2$ by proving the sharp estimate
 $$ \frac{\h^2(\{u>0 \}\cap B)}{\h^2(\{u<0 \}\cap B)} \leq C \log N.$$

 The extra factor of $\log\log \lambda$ is the price payed for using quasiconformal mappings.

\section{Real-analytic Riemannian manifolds and the  work of Donnelly and Fefferman} \label{real-analytic df}
 Donnelly and Fefferman \cite{DF} proved  the Yau conjecture in the case when the metric is real-analytic. 
The work of Donnelly and Fefferman brought many new ideas  to the field. Some of the ideas use the real analyticity of the metric,
some of them work in the smooth setting. In this section we would like to focus on how the real analyticity helps.

\subsection{Real-analyticity in local coordinates.}
   In local coordinates one can think about $\varphi_\lambda$ as of a solution  to the elliptic equation 
\begin{equation} \label{eq:div} 
 \frac{1}{\sqrt {|g|}}\dv(\sqrt {|g|} (g^{ij}) \partial_j \varphi_\lambda) + \lambda \varphi_\lambda=0 
\end{equation}
 in some domain $D$ in $\mathbb{R}^n$. 
  In the  case when the metric is real-analytic
the coefficients of the equation are  real-analytic.

  The following extremely useful idea is due to Donnelly and Fefferman. The idea is using real-analyticity.

 \textbf{Main idea.}
 There is a complex neighborhood $D^*$ of $D$ in $\mathbb{C}^n$, which depends on $g$ and $D$, but does not depend on $\lambda$, such that 
 any solution $\varphi_\lambda$ to \eqref{eq:div} in $D$ has a holomorphic extension  onto $D^*$ with estimate 
\begin{equation} \label{eq:hol ext}
  \sup_{D^*} |\varphi_\lambda| \leq  e^{C\sqrt \lambda} \sup_{D}|\varphi_\lambda|.
\end{equation}

  If we fix a linear elliptic operator $L=div(A\nabla \cdot)$ with real-analitic coefficients in a ball $B\subset \mathbb{R}^n$ with center at the origin, then any solution $u$ of $Lu=0$ is real-analytic \cite{MN57},\cite{J2004},\cite{P39} and moreover the Cauchy estimates for the derivatives of $u$ hold: 
 $$ |D^\alpha u(0)| \leq C \alpha! \sup_{B} |u| / R^{|\alpha|},   $$
  where $R$ depends on the real analytic coefficients of $L$. It implies that $u$ coincides with its Taylor series in some neighborhood of $0$. One can plug complex numbers in the Taylor series of $u$, which naturally defines the holomorhic extension of $u$  in a complex ball $B^*\subset\mathbb{C}^n$ (with smaller radius $R$ than of $B$)  with an estimate: 
   $$\sup_{B^*} u \leq C_L \sup_{B} u.$$
  
However the coefficients of the equation for $\varphi_\lambda$ grow with $\lambda$ and it is  not clear apriori why the domain of holomorphic extension could be chosen independent of $\lambda$.
The original proof \cite{DF} of holomorphic extesnion by Donnelly and Fefferman  was further simplified with the help of the harmonic extension, which allows to pass to an elliptic equation with fixed coefficients. The idea of harmonic extension will be formulated in Section \ref{lift and frequency}. The authors learned this idea from F.-H. Lin.

\subsection{Upper bound in Yau's conjecture  in the real-analytic case.}
The interested reader may start with a simpler case of harmonic functions in $\mathbb{R}^n$ and learn in \cite{Han07} the idea  how to apply the holomorphic extension to upper bounds for nodal sets of harmonic functions.
 The common idea is using the fact that the size of the zero set of holomorphic functions can be estimated in terms of growth of the function.

 \textbf{Complex analysis lemma.}
 Let $f$ be a holomorphic function of one variable in the disc  $\{ |z|<2 \}$
such that $f(0)=1$ and $\sup_{\{ |z|<2 \}}|f| \leq 2^N$ for some number $N$.
 Then 
\begin{equation} \label{eq:complex}
   \textup{ Number of zeroes of $f$ on $ \mathbb{R}\cap \{ |z|<1 \}$ is smaller than $C N$.}
\end{equation}

 The  ideas in the proof of  the upper bound in Yau's conjecture (in the real-analytic case) are the holomorphic extension  \eqref{eq:hol ext} into a complex neighborhood (which does not depend on $\lambda$) and estimate \eqref{eq:complex}
on the number of zeroes of a holomorphic function. 
The whole proof contains technical details. We give a very brief sketch of the technical details.
 First, the doubling index estimate \eqref{eq:di} implies that for any geodesic ball $B$ on the manifold there is a constant $C_1$, which depends
 on the radius of $B$ and on the manifold, such that $$\sup_{B} |\varphi_\lambda| \geq  e^{-C_1\sqrt \lambda}\sup_{M} |\varphi_\lambda|.$$

 So near each point on the manifold one can find a point where the value is not too small.
 Second, consider any point $x$ with value $|\varphi_\lambda(x)|$ at least $e^{-C_1\sqrt \lambda} \sup_{M} |\varphi_\lambda|$. Assume $\sup_{M} |\varphi_\lambda|=1$.
 Looking at $\varphi_\lambda$ in local  coordinates near $x$ one can use the holomorphic extension with estimate and the estimate \eqref{eq:complex} 
to conclude that any segment (in local coordinates) passing through $x$ of length smaller than some constant $c_1$ contains at most $C_2 \sqrt \lambda$ zero points.
The final technical step is to obtain the estimate for $n-1$ dimensional Hausdorff measure of the zero set  from the fact that zero set does not have many intersections with the segments passing through the points, where $|\varphi_\lambda|$ is at least $e^{-C_1\sqrt \lambda}$.

 The idea of the last step is formalized in the next claim.

\textbf{Estimate of Hausdorff measure via intersections with lines.} Fix $n+1$ points $x_1$, $x_2$, \dots $x_{n+1}$  in $\mathbb{R}^n$, such 
 that $x_1$, $x_2$, \dots $x_{n+1}$ do not lie on one $(n-1)$ - dimensional plane.  Suppose that $S$ is a closed set inside of the unit ball $B=\{ |x|<1\}$. If for any line $L$
 passing through at least one point of $x_1$, $x_2$, \dots $x_{n+1}$ the number of points in $L\cap S$ is smaller than a number $N$, then
 $(n-1)$ - dimensional Hausdorff measure of $S$ is smaller than $ CN$, where $C=C(x_1,..., x_{n+1})$ depends on how degenerate the simplex $\{x_1,..., x_{n+1}\}$ is, on the diameter of the ball $B$ (which contains $S$) and on the distance between $B$ and the simplex.
 
\textbf{Remark.}
 If a compact set $S$ in $\mathbb{R}^n$ has a property that any line $L$ passing through  $x_1$ contains at most $1$ point of $S$, then
it is not true that $S$ has a finite $(n-1)$ - dimensional Hausdorff measure. However the last statement becomes true if  the same property
 holds for $n+1$ points  $x_1$, $x_2$, \dots $x_{n+1}$, which do not lie on the same hyperplane.
 
\textbf{Estimate for the number of balls of size $\frac{1}{\sqrt \lambda}$ with large doubling index.} 
One can cover  $M$ by $\sim \lambda^{n/2}$ balls $B_j$ of radius $ \frac{1}{\sqrt \lambda} $ in such a way
that  every point of $M$ is covered at least once and at most $C$ times. 

Donnelly and Fefferman proved that for at least half of $B_j$ the doubling index of $\varphi_\lambda$ in $B_j$ is bounded by some constant $C_1$, which does not depend on $\lambda$. 
One can also replace the word ``half'' by $99/100$ and the statement above will remain true, but $C_1$ will become larger.
In other words for most of balls of size $ \frac{1}{\sqrt \lambda} $ the doubling index of $\varphi_\lambda$ is controlled.

\textbf{Remark.}
The latter statement has beautiful corollaries: lower bound in Yau's conjecture and the quasi-symmetry conjecture.

 \subsection{Lower bound in Yau's conjecture  in the real-analytic case.}

 The zero set of $\varphi_\lambda$ is $\frac{C}{\sqrt \lambda}$ dense on $M$. One can cover
 $n$-dimensional closed manifold $M$ by $\sim \lambda^{n/2}$ balls $B_j$ of radius $\sim \frac{1}{\sqrt \lambda} $ in such a way
that $\varphi_\lambda$ is zero at the centers of the balls $B_j$ and every point of $M$ is covered less than $C$ times.

 Donnelly and Fefferman showed that at least half of $B_j$ have a doubling index smaller than $C_1$. Furthermore they
 showed that if a ball $B_j$ of radius $\frac{C}{\sqrt \lambda}$ has a doubling index smaller than $C_1$ and $\varphi_\lambda$ is zero at the center of $B_j$, then
 $$\h^{n-1}(B_j \cap \{\varphi_\lambda=0\}) \geq \frac{c_1}{(\sqrt \lambda)^{n-1}} $$
and 
  $$C_2 > \frac{\h^{n}(B_j \cap \{\varphi_\lambda > 0\})}{\h^{n}(B_j \cap \{\varphi_\lambda < 0\})} \geq c_2>0.$$

 Since the total number of such $B_j$ is comparable to $(\sqrt \lambda)^{n}$, it yields the lower bound in Yau's conjecture
 and the quasi-symmetry conjecture in the real-analytic case.

 \textbf{Remark on the lower bound in the real-analytic case.} 
 The proof of the lower bound in Yau's conjecture is more elaborate than of the upper bound.
 The most interesting part is to how show that for at least half of balls $B_j$ of size $ \frac{C}{\sqrt \lambda} $ the doubling index of $\varphi_\lambda$ is bounded. 
 The proof of the latter fact is also using the holomorhpic extension \eqref{eq:hol ext} with the growth estimate $e^{C\sqrt \lambda}$ . We do not dare to explain the complete plan of the proof, but
 we would like to mention one useful statement, which helps to control oscillations of holomorphic functions in terms of growth. It gives a hint why we should believe that  the doubling index for half of the balls is bounded.

 \textbf{Complex/Harmonic analysis lemma}( Proposition 5.1 in \cite{DF}).
 Fix $\varepsilon >0$.
Let $f$ be a holomorphic function of one complex variable in the disc  $\{ |z|<2 \}$ 
such that $f(x)$ is real   for $x \in [-2,2]$, $f(0)=1$ and $\sup_{\{ |z|<2 \}}|f| \leq 2^N$ for some integer number $N>1$. Split the interval $[-1,1]$ into $N$ equal intervals $Q_\nu$ of length $\frac{1}{N}$.
 Then there is a set $E \subset [-1,1]$ of measure less than $\varepsilon$ such that 
 $$ |\log f^2(x) - \log (\frac{1}{|Q_\nu|} \int_{Q_\nu} f^2)| \leq C_{\varepsilon} $$
for any $x\in Q_\nu \setminus E$. Here $C_{\varepsilon}$ depends only on $\varepsilon $.

\textbf{Remark.} The statement above allows to control oscillations of holomorphic functions in terms of growth.
 It suggests that if a holomorphic function grows slower than a polynomial of degree $N$, then 
it behaves nice on most of the intervals of length $1/N$.
In particular if $\varepsilon$ is small enough, then the lemma used twice for $N$ and $2N$  gives
$$\frac{1}{|Q_\nu|} \int_{Q_\nu} f^2 \leq C \frac{1}{|\frac{1}{2}Q_\nu|} \int_{\frac{1}{2}Q_\nu} f^2 $$
for at least half of $Q_\nu$. In other words $L^2$ doubling index for $f$ is bounded for a big portion of intervals of length $1/2N$.

\section{Norm estimates of eigenfunctions and their applications for the lower bound in the Yau conjecture}
\subsection{Weyl's law}  A classical result on the spectrum of the Laplace operator $\Delta$ on a compact Riemannian manifold $M$ is the Weyl asymptotic law.  Let $N(\lambda',\lambda'')$ denote the number of eigenvalues $\mu$ of the operator of $-\Delta$ such that $\lambda'\le\mu<\lambda''$. Then
\[N(0,\lambda)=c_n\lambda^{n/2} vol(M)+O(\lambda^{(n-1)/2}),\]
where $c_n$ is the constant that depends only on the dimension. It implies that $N((k-1)^2,k^2)$ is comparable to $k^{n-1}$. 

Motivated by the study of projections in spherical harmonics, Sogge \cite{S88} considered projections $P_k$ of $L^2(M)$ onto subspaces generated by eigenfunctions with eigenvalues $\mu\in[(k-1)^2,k^2)$. He obtained sharp inequalities of the form
\[\|P_kf\|_{q}\le k^{\sigma(p,q,n)}\|f\|_{p},\]
where $1\le p\le 2$ and $q=2$ or $p'=p/(p-1)$. These inequalities imply $L^p$ norm estimates for eigenfunctions. In particular, if    $\varphi_\lambda$ is an eigenfunction  and  $p=2(n+1)/(n-1)$ then
\begin{equation}\label{eq:Lp}
\|\varphi_\lambda\|_p\le C\lambda^{1/(2p)}\|\varphi_\lambda\|_2.
\end{equation}

\subsection{Lower estimate of the zero set on the smooth case by Colding and Minicozzi} Inspired by the ideas of Donnelly and Fefferman, Colding and Minicozzi \cite{CMII} proved a lower bound  for the size of the nodal set by finding many balls of size $C/\sqrt{\lambda}$ with bounded doubling index. They covered  the manifold $M$  by balls of radius $C/\sqrt{\lambda}$ in such a way that $\varphi_\lambda$ is zero at the center of each ball and each point is covered by not more than $C_1$ balls. As in the real analytic case it is also true in the smooth case that one can control the size of the zero set in balls, where the doubling index is smaller than a fixed numerical  constant $D$. Let's call such balls good and denote by $\mathcal{B}$ the collection of good balls.
In each good ball $B$ the eigenfunction $\varphi_\lambda$
cannot oscillate too fast since we control its growth properties, and  one can inscribe a ball 
of radius $c_D/\sqrt \lambda$ in $B \cap \{\varphi_\lambda >0 \}$  and a ball of radius $c_D/\sqrt \lambda$ in $B \cap \{\varphi_\lambda <0 \}$. Every segment connecting these balls has a sign change of $\varphi_\lambda$ and that therefore
$$\h^{n-1}(Z(\varphi_\lambda)\cap B)\ge c\lambda^{-(n-1)/2}.$$

 The elegant idea \cite{CMII} states that most of the $L^2$-mass of the function is concentrated in good balls:
\[\sum_{B\in\mathcal{B}}\int_B|\varphi_\lambda|^2\ge c_1\|\varphi_\lambda\|_2^2.\] 

 This is true because the sum  over bad balls satisfies
\[\sum_{B\notin\mathcal{B}}\int_B|\varphi_\lambda|^2\le \frac{C}{D} \|\varphi_\lambda\|_2^2.\]

Let $G=\cup_{B\in\mathcal{B}} B$.
Since each point is covered at most $C_1$ times we have 
\[\int_{G}|\varphi_\lambda|^2\ge c\|\varphi_\lambda\|_2^2.\] 
 Then the H\"older inequality gives
\[c\|\varphi_\lambda\|_2^2\le\int_{G}|\varphi_\lambda|^2\le \left(\int_{M}|\varphi_\lambda|^{2(n+1)/(n-1)}\right)^{(n-1)/(n+1)}|G|^{2/(n+1)}.\]
To estimate the number of good balls, Colding and Minicozzi applied the $L^p$ bound  \eqref{eq:Lp} of Sogge, it implies $|G|\ge c_2\lambda^{-(n-1)/4}$.  $G$ is a union of good balls of radius $C/\sqrt \lambda$. Hence the number of the good balls is at least $c_3\lambda^{(n+1)/4}$. This leads to the following lower bound of the zero set
\[\h^{n-1}(Z(\varphi_\lambda))\ge c_4\lambda^{-(n-3)/4}.\]

\subsection{Lower estimate of the zero set by Sogge and Zelditch} Another approach to the estimate of the size of the zero set from below was suggested \cite{SZ11},\cite{SZ12} by Sogge and Zelditch.
Their starting point is the corollary of the Green formula applied to nodal domains 
\[\lambda\int_M|\varphi_\lambda|=2\int_{Z(\varphi_{\lambda})}|\nabla\varphi_\lambda|.\]
It immediately implies that
\[\h^{n-1}(Z(\varphi_\lambda))\ge \lambda \frac{\|\varphi_\lambda\|_1}{ \|\nabla\varphi_\lambda\|_\infty}. \]
Rescaling the equation $\Delta\varphi_\lambda+\lambda\varphi_\lambda=0$ from balls of radius $c/\sqrt{\lambda}$ to unit balls and applying standard elliptic estimates it is not difficult to prove that
$\|\nabla\varphi_\lambda\|_\infty\le C\sqrt{\lambda}\|\varphi_\lambda\|_\infty.$ Finally  the inequality 
\[\|\varphi_\lambda\|_\infty\le C\lambda^{(n-1)/4}\|\varphi_\lambda\|_1\]  implies the estimate $\h^{n-1}(Z(\varphi_\lambda))\ge c\lambda^{-(n-3)/4}.$ The inequality between the $L^1$ and $L^\infty$ norms of eigenfunctions is non-trivial.

\subsection*{Remark} We would like to mention that the third proof of the lower bound \cite{CMII},\cite{SZ12} (mentioned in the previous subsection) was given \cite{St14} by Steinerberger, who applied the heat flow to the eigenfunctions. His approach is  using $L^\infty$ bounds for eigenfunctions. The heat kernel approach in \cite{St14} was further generalized to linear combinations of eigenfunctions  \cite{St17} and it shows that if $f$ is a linear combination of eigenfunctions with eigenvalues in the interval $(\lambda, 10 \lambda)$, then either the $L^\infty$ norm of $f$ is large compared to $L^1$ norm or $f$ has a large zero set.

 Other works on the lower bounds include \cite{M11,HS12}.

 \section{From eigenfunctions to solutions of elliptic PDEs} \label{lift and frequency}
\subsection{Harmonic extension of eigenfunctions} Some of the questions on the Laplace eigenfunctions can be reduced to questions on harmonic functions on Riemannian manifolds. If $\varphi_\lambda$ is an eigenfunction on $M$ satisfying $\Delta \varphi_\lambda+\lambda\varphi_\lambda=0$ then the function $$u(x,t)=\varphi_\lambda(x)e^{\sqrt{\lambda}t}$$ is defined on $M\times\R$ and is  harmonic  with respect to the natural metric on the product manifold. 
 The zero set of the new function is the cylinder over the zero set of the old function:
 $$Z_u=Z_{\varphi_\lambda} \times \mathbb{R}.$$
The equation for $u$ is now independent of $\lambda$, the growth of $u$ in the variable $t$ carries the information on the eigenvalue. The idea of harmonic extension allows
 to simplify several steps in the proof of the upper bound in \cite{DF}, including the holomorphic extension with estimate and the proof of the doubling index estimate $C\sqrt \lambda$.  

\subsection{The frequency function}\label{ss:freq} The doubling index of a solution of elliptic PDE is connected to the so-called frequency function. We start with a harmonic function $h$ in a subdomain $\Omega$ of the Euclidean space and for each ball $B=B_r(x)\subset\Omega$ define the following quantities
\[H_h(x,r)=\frac{1}{|\partial B|}\int_{\partial B}h^2, \quad \F_h(x,r)=r\frac{d}{dr}\log H(x,r).\] 
It was known to Agmon \cite{Ag} and Almgen \cite{Al} that $\F_h(x,r)$ is an increasing function of $r$ and thus the function $t\to \log H_h(x,e^t)$ is convex; $\F_h$ is called the frequency function of $h$. Garofalo and Lin \cite{GL}
 showed that a similar almost monotonicity inequality holds for solutions of second order elliptic PDEs in divergence form with Lipschitz coefficients, which has many applications to nodal sets on smooth manifolds.  We omit the accurate definition of the frequency in that setting, but we would like to describe the relation of this term to the doubling index defined in \ref{ss:di}.
The doubling index and frequency are almost synonyms:
one of them deals with the growth of $L^2$ norm another one with $L^\infty$ norm.   The definition of the frequency and the monotonicity property lead to the following inequalities 
\[
\F_h(x,r)\le\log_2\frac{H_h(x,2r)}{H_h(x,r)}\le \F_h(x,2r)
\]
for any ball $B=B(x,r)$ such that $B(x,2r)\subset\Omega$.
At the same time the standard elliptic estimates imply that one can compare $L^2$ and $L^\infty$ norms of harmonic functions if we are allowed to enlarge the radius of the ball:
\[H_h(x,r)\le\max_{B_r(x)}|h|\le CH_h(x,3r/2).\] 
This explains the connection between the frequency and the doubling index of solutions of second order elliptic equations,
\[c\F_h(x,r/3)\le N_h(B_r(x))\le C\F_h(x,3r).\]  

\subsection{Applications to eigenfunctions} We consider the lift $u(x,t)$ of an eigenfunction $\varphi_\lambda$  as  a solution of the second order elliptic PDE in divergence form. The dependence of $u$ on $t$ encodes the growth properties: 
$$ \max_{M\times[-2r,2r]} |u| = e^{\sqrt \lambda r} \max_{M\times[-r,r]} |u|.$$
It is not difficult to deduce from the almost monotonicity property of the frequency function or of the doubling index that the frequency and the doubling index  of $u$ in any ball in $M\times [-1,1]$ are bounded by $C\sqrt{\lambda}$. This approach  gives a  new proof of the doubling index estimate $C\sqrt \lambda$ of eigenfunctions.
For harmonic functions on Riemannian manifolds the frequency controls the order of vanishing and also the size of the zero set, see \cite{Han07},\cite{HS},\cite{L1}.

\section{Propagation of smallness}
Harmonic functions on Riemannian manifold share some properties of holomorphic functions when the metric is smooth but not necessarily real-analytic. An important illustration is the (weak) unique continuation principle for such functions proved by Cordes \cite{Co} and Aronszajn \cite{Ar}. If a harmonic function defined on some domain $\Omega$ is zero on an open subset of $\Omega$ then it is zero on the whole domain. The statement is true for solutions of elliptic equations in divergence form with Lipschitz coefficients,  but it fails \cite{Mi73} for Holder continuous coefficients in dimension larger than $2$. We need two quantitative versions of this principle which we refer to as propagation of smallness. 

 \subsection{Three ball theorem} The first quantitative version of the statement is the Hadamard three-circle theorem for holomorphic functions on the complex plane. Let $f$ be holomorphic in some neighborhood of the origin, define $M(r)=\max_{|z|=r}|f(z)|$, then
\[M(r_1)\le M(r_0)^\alpha M(r_2)^{1-\alpha},\quad{\text{where}}\ r_0<r_1<r_2,\ r_1=r_0^\alpha r_2^{1-\alpha}.\]
For a function $h$ harmonic in Euclidean metric the same statement becomes true \cite{Han07},\cite{Ag},\cite{Al} if one replaces $M(r)$ by the $L^2$-average $H_h(x,r)=
|\partial B_r(x)|^{-1}\int_{\partial B_r(x)}|h|^2$. It is equivalent to  the convexity property discussed in section \ref{ss:freq}. Further, comparison of $L^2$ and $L^\infty$ norms yields the following inequality
\begin{equation}\label{eq:3B}
\sup_{|x|\le r_1}|h(x)|\le C\sup_{|x|\le r_0}|h(x)|^\alpha\sup_{|x|\le r_2}|h(x)|^{1-\alpha},\quad \alpha=\alpha(r_0/r_1,r_2/r_1)\in(0,1).
\end{equation}
This result was generalized to solutions of elliptic equations with non-analytic coefficients by Landis \cite{L63}. It holds for solution $u$ to elliptic equations in divergence form $\dv(A\nabla u)=0$ with Lipschitz coefficients.

\subsection{Question of Landis} \label{Landis spheres} Landis asked whether the inequality \eqref{eq:3B} remains true when the smallest ball $\{|x|\le r_0\}$ is replaced by a wild set of positive measure. Suppose that $u$ is a solution of the equation $\dv(A\nabla u)=0$ in some ball $2B$ such that $|u|$ is bounded by one in $2B$ and let $E$ be a measurable subset of $B$ with positive measure. Assume also that $u$ is small on $E$, $\sup_E|u|\le\varepsilon$. The question of Landis is whether the inequality
\[\sup_B|u|\le C\varepsilon^\alpha\]
holds with some $C>0$ and $\alpha\in(0,1)$ that depend on the equation and on the volume $|E|>0$, but not on $u$ or the geometry of $E$. For the case of real analytic coefficients the affirmative answer was given by Nadirashvili \cite{N79}. For the case of smooth coefficients partial advances towards the question of Landis were obtained by Nadirashvili \cite{N86} and by Vessella \cite{V00}. In the real-analytic case it can be obtained with the help of holomorphic extension with estimate. For any holomorphic function $f$ (of one complex variable)  $\log|f|$ is subharmonic, which is a powerful unique continuation property. It allows to get the affirmative answer to the question of Landis problem in the real-analytic case.

 The positive answer to the question of Landis was recently obtained for elliptic equations of the form $\dv(A\nabla u)=0$ with Lipschitz coefficients. The proof is based on a simple version of multiscale iteration, see the lecture notes  \cite{LM2} for the mini-course at PCMI 2018 and also \cite{LM1} for further discussion of the classical problems on unique continuation. The recent proof  also yields the bound for the BMO-norm of $\log|u|$ in terms of the doubling index of $u$.  As the corollary one can obtain the estimate of the BMO-norm of $\log|\varphi_\lambda|$  for eigenfunctions that was discussed in Section \ref{ss:BMO}.

\subsection{Quantitative Cauchy uniqueness theorem} 
The Cauchy uniqueness property for second order elliptic PDEs with Lipschitz coefficients states that if $\dv(A\nabla u)=0$ in some domain $\Omega$, $u\in C^1(\bar{\Omega})$ and the solution $u$ and its normal derivative $u_n$ vanish on a relatively open part $\Gamma$ of $\partial\Omega$ then $u=0$ in $\Omega$.  We are looking for a quantitative version of this statement, which is called conditional stability of the Cauchy problem. For the history of the question that goes back to Hadamard we refer the reader to the survey \cite{ARRV}. 

We formulate the quantitative result from \cite{ARRV} in simple geometric situation. Suppose that $u$ is a solution of an elliptic equation $\dv(A\nabla u)=0$ in the half-ball $B_{2r}^+=\{x=(x_1,...,x_n): |x|<2r, x_n>0\}$ and is $C^1$ smooth up to the boundary. Let $\Gamma$ be the flat part of the boundary of $B_{2r}^+$, $\Gamma=\{x:|x|<r, x_n=0\}$. We assume   $\sup_{B_{2r}^+}|u|\le 1$. If the Cauchy data of $u$ is small on $\Gamma$, then the solution is small on each smaller half-ball. Namely, if $\sup_\Gamma(|u|+|\nabla u|)\le\varepsilon$ then
\[\sup_{B_r^+}|u|\le C\varepsilon^{\alpha},
\]
for some $C>0, \alpha\in(0,1) $ that depend on the coefficients of the equation and on $r$ but not on $u$.

\subsection{Old open question on the Cauchy uniqueness problem.} \label{Cauchy problem} Here we mention the second question that Nadirashvili suggested \cite{N97} to focus on to solve the Yau conjecture. 

Assume $u$ is a harmonic function in the unit ball $B\subset \mathbb{R}^3$ and $u$ is $C^\infty$-smooth in the closed ball $\overline{B}$.
 Let  $S \subset \partial B $ be any closed set with positive area.
Is it true that $\nabla u = 0$ on $S$ implies $\nabla u \equiv 0$?

 For the less smooth  class of functions  $C^{1+\varepsilon}(\overline{B_1} )$ there is a striking counterexample \cite{BW}, \cite{W}. The attempts to construct $C^2$-smooth counterexamples were not successful.

\section{Zeroes and singular sets of solutions of elliptic PDEs with smooth coefficients}
\subsection{Structure of the zero set and the estimate of Hardt and Simon} 
Recall that the singular set $S(\varphi)$ of a function $\varphi$ is the set $\{x: \varphi(x)=0, \nabla \varphi(x)=0 \}$ and the critical set is defined by $\{ x :\nabla \varphi(x)=0 \}$.
The structure of the zero sets and critical sets of eigenfunctions and more general solutions of elliptic PDEs is a delicate subject especially when the metric or the coefficients of the equations are not real-analytic. 

  Bers \cite{B55} studied the local behavior of solutions of linear elliptic PDE by looking at the first term in the Taylor expansion of the solution. In particular the work of Bers implies
that for Laplace eigenfuncitons on surfaces the nodal set
is a union of curves with equiangular intersections.
 Caffarelli and Friedman \cite{CF}   showed that in dimension $n$ the singular sets of some linear and semilinear elliptic equations have Hausdorff dimension at most $n-2$.  Hardt and Simon \cite{HS} proved upper bounds for Hausdorff measure of nodal sets in terms of growth.  Their results imply that if $u$ is a solution to a second order elliptic equation in divergence form with Lipschitz coefficients  in a fixed ball $2B$  then
\[\h^{n-1}(Z(u)\cap B)\le CN^{CN},\]
where $$N=N_u(B)=\log \frac{\sup_{2B}|u|}{\sup_{B}|u|}$$ is the doubling index.

\subsection{Harmonic counterpart of the estimate from above in Yau's conjecture}
The estimate of Hardt and Simon was recently improved in \cite{L1}, where it was shown that for solutions of second order elliptic equations in divergence form $\dv (A \nabla u) =0$ with smooth coefficients there exists $a=a(n)$ such that  
\[\h^{n-1}(Z(u)\cap B)\le CN^{a(n)},\quad N=N_u(B),\]
where $C$ depends on the coefficients of the equation in $2B$, but not on $u$.  The polynomial upper bound in the Yau conjecture follows from this inequality and the lifting trick. 

The  upper bound in Yau's conjecture \[\h^{n-1}(Z(\varphi_\lambda))\le C\sqrt{\lambda}\] would follow from the linear estimate on the size of the zero set in terms of the doubling index: $$\h^{n-1}(Z(u)\cap B)\le CN.$$ 

\subsection{Estimates of the singular set} The quantitative estimates of the singular and critical sets of solutions of elliptic PDEs is an interesting and developing topic. For harmonic funcitons the critical set is of codimension two and has locally finite $n-2$-dimensional  Hausdorff measure.  The proofs \cite{HHHN},\cite{HHN96},\cite{HHN97} show that the measure of the singular set can be estimated in terms of the doubling index for the gradient. The explicit bound
\[\h^{n-2}(S(u)\cap B)\le C^{N^2},\quad N=N_u(B)\]
 was recently obtained \cite{NV} by  Naber and Valtorta.
The method to study stratification properties of the singular sets involves the notion of the effective singular set and is also explained in \cite{CNV} by Cheeger, Naber and Valtorta. 

 Examples of singular sets of harmonic functions in dimension three suggest much stronger bound $\h^{n-2}(S(u)\cap B)\le C N^2$, which was conjectured by F.-H. Lin.  It is not known whether this estimate holds even for harmonic functions in $\R^3$.

\section{On the proof of the polynomial upper bound} \label{upper bound}
 Let $Q_0$ be a cube in $\mathbb{R}^n$. Consider a solution $u$ to an elliptic equation in divergence form $\dv(A\nabla u)=0$ with Lipschitz coefficients $A$ 
in a cube $3 Q_0$. By $N_u(Q_0)$ we denote the doubling index of $u$ in a cube $Q_0$:
$$N_u(Q_0)= \log_2 \frac{\sup_{2Q_0}|u|}{\sup_{Q_0}|u|}.$$
 The proof \cite{L1} of the polynomial upper bound
\begin{equation} 
\h^{n-1}(Z_u \cap Q_0) \leq CN^{C_n}_u(Q_0).
\end{equation}
 is a multiscale iteration argument in its essence. 
Typically such arguments are cumbersome, but 
 we hide all iterations in one notation ($F(N)$) so that the reader does not see  iterations at all.

Assume that any subcube  of $Q_0$ has a doubling index not greater than a number $N$ and 
we would like to find the smallest upper bound $F(N)$ such that 
\begin{equation} \label{eq:upper}
\h^{n-1}(Z_u \cap Q) \leq F(N) s^{n-1}(Q)
\end{equation}
 for any subcube $Q$ of $Q_0$,
 where $s(Q)$ is the side length of $Q$.
  Hardt and Simon proved  in \eqref{eq:upper}  that $F(N)\leq CN^{CN}$, in particular $F(N)$ is finite. 
 While our goal is to show that  $F(N) \leq CN^{C}$.

 Assume that a cube $Q$ has side length $1$. If $N$ is small we can use the bound by Hardt and Simon. When $N$ is  large we can chop the unit cube $Q$ into $K^n$ equal  smaller cubes $q_i$ of size $1/K$,
  and  we may hope that many of $q_i$ will have a doubling index smaller than the doubling index of $Q$ and in fact it happens and helps.

\textbf{Lemma on the distribution of doubling index.} If $K$ and $N$ are sufficiently large, then there are at least  $K^n- \frac{1}{2}K^{n-1}$
good cubes $q_i$ such that $N(q_i) \leq N/2$ and for any subcube $\tilde q$ of good $q_i$ the doubling index of $\tilde q $ is also smaller than $N/2$.

Consider the size of the zero set inside of each cube $q_i$. For good cubes we have a good bound   $F(N/2)\frac{1}{K^{n-1}}$ and for the cubes, which are not good, we only have a bound $F(N)\frac{1}{K^{n-1}}$.
 But the number of bad cubes (which are not good)  is smaller than $\frac{1}{2}K^{n-1}$.
 The latter lemma implies the recursive inequality:
$$F(N) \leq 2K F(N/2),$$
which yields the polynomial upper bound.

\textbf{Remark}. For the applications of the lemma on the distribution of doubling index it is important that the constant $\frac{1}{2}$ in lemma is smaller 
than $1$. If there was  $1$ instead of $\frac{1}{2}$ we would not be able to obtain any upper bound on Hausdorff measure.

The proof of lemma on the distribution of doubling indices is using techniques of quantitative unique continuation. 
There are two main statement in the proof: the simplex lemma and the hyperplane lemma. For the sake of simplicity we formulate them for ordinary harmonic functions in $\R^3$.

\textbf{Simplex lemma}. Suppose  4 points $x_1,x_2,x_3, x_4$ in $\R^3$ form a non-degenerate simplex $S$ with sides at least $1$. 
Define the width of $S$ as
 the minimal distance  between  pairs of parallel planes in $\R^3$ such that  $S$ is  between two planes.  Let $a>0$ and assume that $\frac{\textup{width}(S)}{\textup{diam}(S)} \geq a$. 
There exist positive constants $c=c(a),C=C(a),k= k(a)$ such that the following holds.  
 If $u$ is a harmonic function in $\R^3$ such that $N_u(B_1(x_i)) \geq N$ for $i=1,2,3,4$, 
then for the barycenter $x_0$ of $S$ the doubling index of $u$ in $B_{k \textup{ diam}(S)}(x_0)$ is at least   $N(1+c)-C$.

\textbf{Hyperplane lemma}. Let $u$ be a harmonic function in $\R^3$, let $A>100$ be  an integer and $N\geq 2$. 
Consider a finite lattice of points $L_A:\{ (i,j,0): i=- A, \dots A, j=-A, \dots A \}$. If $N$ and $A$ are sufficiently large and $B_1(x) \geq N$ for each $x\in L_A$, then $N(B_A(0))\geq 2N - C$.

 The proof of the simplex lemma is using the monotonicity property of the frequency and the proof of the hyperplane lemma relies on the quantitative Cauchy uniqueness property.

\section{Lower bound in Yau's conjecture} \label{lower bound}
\subsection{Reduction to Nadirashvili's conjecture}
  The proof of the lower bound in the Yau conjecture is using the fact that the zero set of $\varphi_\lambda$ is $C/\sqrt \lambda$ dense.
 The manifold $M$ can be covered by $\sim \lambda^{n/2}$ balls $B_j$ of radius $\sim \frac{1}{\sqrt \lambda} $ in such a way
that $\varphi_\lambda$ is zero at the centers of the balls $B_j$ and every point of $M$ is covered less than $C$ times. 
Donnelly and Fefferman proved that in the real analytic case at least half of $B_j$ have a bounded doubling index, but we don't know 
whether it holds in the smooth case and we follow another way suggested by Nadirashvili. We prove that in each ball $B_j$ (with zero at its center) the zero set 
of $\varphi_\lambda$ has $(n-1)$-dimensional Hausdorff measure at least $c\lambda^{-(n-1)/2}$. Since each ball is covered at most $C$ times and the total number of balls is comparable to $\lambda^{n/2}$ the lower bound in Yau's conjecture  follows.

 Nadirashvili proposed a conjecture about harmonic functions in order to attack the lower bound in Yau's conjecture.

{\bf Conjecture of Nadirashvili} (proved in \cite{L2}). There exists a constant $c>0$ such that for any harmonic function $u$  in a unit ball $B$ in $\mathbb{R}^3$ with zero at the center of $B$
 the area of zero set of $u$ in $B$ is larger than $c$.

 In dimension two a similar question is not difficult. Zero set of any non-zero harmonic function in a unit disc is a union of analytic curves and due to the maximum principle each zero curve is not allowed to have loops and therefore there is a zero curve, which connects the center of of the disc with the boundary of the disc and has a length bigger than the radius of the disc. 

In dimensions three and higher the conjecture is true, but the recent proof given in \cite{L2} is complicated and works for solutions of more general elliptic equations. The proof does not use real-analyticity. No simple  proof is known for the case of harmonic functions in $\R^3$.   

To prove the lower bound in Yau's conjecture we need the rescaled version of Nadirashvili's conjecture for elliptic equations.

{\bf Theorem}.
If $u$ is  a solution of a second order uniformly elliptic equation $\dv(A\nabla u)=0$ with Lipschitz coefficients  in the unit ball $B_1(0)\subset\R^n$ and $u(0)=0$ then \[\h^{n-1}(Z(u)\cap B_r(0))\ge c r^{n-1}, r\in (0,1),\] where $c$ depends on the  coefficients of elliptic equation but does not depend on the solution $u$.

Combining the last theorem with the lifting trick one can show that  in each ball $B_j$ (with zero of $\varphi_\lambda$ at its center and of radius $\sim\frac{1}{\sqrt \lambda}$) the zero set 
of $\varphi_\lambda$ has $(n-1)$-dimensional Hausdorff measure at least $c\lambda^{-(n-1)/2}$.  

 \subsection{On the proof of Nadirashvili's conjecture.}
The proof is not simple and quite long. We would like to explain in this section the logic of the proof and the key points.
We took a liberty to add the notion of stable growth to the proof, which simplifies understanding. In the original text \cite{L2} the words "stable growth" were not used.

Everywhere in this section $u$ is assumed to be a harmonic function in $\R^3$.
We would like to start with a very non-sharp claim.

\textbf{Claim 1}. Let $B_1$ be a unit ball in $\R^3$. If $u(0)=0$ and $N_u(B_1) \leq N$, then 
$$ \h^2(Z_u \cap B_1) \geq \frac{c}{N^{2}}.$$

 The claim is not difficult to prove. It is true that one can inscribe a ball $b_+$ of radius $\frac{c_1}{N}$ in $\{u>0\} \cap B_1$
 and a ball $b_-$ of radius $\frac{c_1}{N}$ in $\{u<0\}\cap B_1$ (see \cite{LM} for details). Every segment connecting $b_+$ and  $b_-$
 has a zero point. That implies $ H^2(Z_u \cap B_1) \geq \frac{c}{N^{2}}.$
 
 The bound  above becomes worse as $N \to \infty$.  Our goal is to obtain a uniform lower bound, which does not depend on $N$.

We will say that  $u$ has a stable growth  in a ball $B$ if 
$$ N_u(B) \leq 1000 N_u(\frac{1}{4}B).$$
The number $1000$ is a fixed sufficiently large constant.

 We will say that  $u$ has a stable growth of order $N$  in a ball $B$ 
if  $N_u(\frac{1}{4}B) \geq N$ and $N_u(B) \leq 1000 N$. 

\textbf{Remark}. Note that by almost monotinicty of the doubling index the  opposite estimate holds $$  N_u(\frac{1}{4}B) \leq C N_u(B).$$
 Here is the key lemma used in the proof of uniform lower bound.

\textbf{Key lemma} (simplified version of Proposition 6.1 in \cite{L2}).
There is a sufficiently large number $N_0$  such that the following holds.
Let $B=B_r(x)$ be a ball of radius $r$ in $\mathbb{R}^3$.  If $N>N_0$ and a harmonic function $u$ has a stable growth of order $N$  in $B$, then 
there exist $c [\sqrt N]^{n-1} 2^{c \log N / \log \log N }$ disjoint balls $B_j$ in $B$ of radius $r/\sqrt N$ such that $u$ is zero 
at the centers of $B_j$.

\textbf{Remark.} 
 The fact that in the key lemma  the constant $2^{c \log N / \log \log N }$ is larger than 1  gives us a hint (but not immediate proof) that
 the bigger the doubling index, the better lower bounds should be. In fact it is true that the bigger $N_u(\frac{1}{4}B)$, the bigger $\h^{2}(Z_u \cap B)$ should be.  

 Just like the proof of the polynomial upper bound, the proof of Nadirashvili's conjecture is also a multiscale  argument in its nature, and again we will hide all multiscale iterations in one notation.
We define $$F(N)= \inf \limits_* \frac {\h^2(Z_u \cap B_r(x))}{r^2},$$
where the infimum is taken over all harmonic functions in $\mathbb{R}^3$ and
 over all balls  $B_r(x)$ such that
\begin{enumerate}
\item[(i)] for any ball $b \subset B_1$ the doubling index of $u$ in $b$ is not greater than $N$,
\item[(ii)] $B_r(x) \subset B_1$,
\item[(iii)] $u(x)=0$.
\end{enumerate}
The rescaled version of Claim 1  gives an estimate $F(N) \geq c/N^2$ and our goal is to show that $F(N) \geq c>0$.

\textbf{Logic of the proof of uniform lower bound using the key lemma.}
Fix $N$ and consider $u$ and a ball $B=B_r(x)$  such that the conditions (i),(ii),(iii) hold and $F(N)$ is almost achieved on $u$:
$$ \frac{\h^2(Z_u \cap B)}{r^2} \leq 2 F(N).$$
  If $N_u(\frac{1}{4}B) \leq N_0$ , we can use Claim 1 to conclude $$F(N) \geq  \frac{\h^2(Z_u \cap B)}{2 r^2} \geq c_2.$$
  When $N_u(\frac{1}{4}B)$ is sufficiently large we  will show that $\frac{\h^2(Z_u \cap B)}{2 r^2} > 2F(N)$ and therefore will arrive to contradiction.

   If  $u$ has a stable growth in $B$, namely $ N_u(B) \leq 1000 N_u(\frac{1}{4}B)$, then we could denote $N_u(\frac{1}{4}B)$ by $\tilde N$ and the  key lemma would imply  that there exist $c [\sqrt {\tilde N}]^{n-1} 2^{c \log \tilde N / \log \log \tilde N }$  disjoint balls $B_j$ in $B$ of radius $r/ \sqrt {\tilde N}$ such that $u$ is zero  at the centers of $B_j$. For each of $B_j$ we know $$\h^2(B_j \cap Z_u) \geq F(N) \left( \frac{r}{\sqrt {\tilde N}} \right)^2.$$
 Since the balls $B_j$ are disjoint   and in $B$ we get 
 $$ \frac{\h^2(B \cap Z_u)}{r^2} \geq \sum \frac{\h^2(B_j \cap Z_u)}{r^2}  \geq c F(N)  2^{c \log \tilde N / \log \log \tilde N }  > 2F(N).$$
The contradiction is obtained. 
  However there is an obstacle to directly apply the key lemma
 because  it is not necessarily true that $u$ has a stable growth in $B$.

 But there is a smaller ball in $B$ with stable growth of $u$. 

\textbf{Lemma on stable growth (Follows from Lemmas 4.1,4.2 in \cite{L2})}. If $N_u(\frac{1}{4}B)$ is sufficently large, then there is  a ball $\tilde B \subset B$ and a number $\tilde N \geq c_3 N_u(\frac{1}{4}B) $
such the radii of $B$ and  $\tilde B$ are related by
$$r(\tilde B)= \frac{c_4 r(B)}{\log^2 \tilde N}$$  and  $u$ has a stable growth of order  $\frac{\tilde N}{\log^2 \tilde N}$ in $\tilde B$.

\textbf{Remark.} The proof of the lemma on stable growth is using the monotonicity property of the frequency and the following fact on monotonic functions.

 \textbf{Fact}. If $\beta(r)$ is any increasing function on $[0,1]$ with $\beta(0) >2$, then 
 there is a number $N\geq 2$ and an interval $I \subset [0,1] $  of length $\frac{c}{\log^2 N}$
 such that $$N \leq \beta(r) \leq 2N \textup{ for } r \in I.$$

  Combining the lemma on stable growth and the key lemma we come to the same conclusion that 
$$ \frac{\h^2(B \cap Z_u)}{r^2} \geq c_5 F(N)  2^{c \log \tilde N / \log \log \tilde N }/ \log^4 \tilde N > 2F(N) .$$

 We finished the attempt to explain how   the key lemma implies Nadirashvili's conjecture. 
 
\textbf{ On the proof of the key lemma.}
In this section we present a plan of the proof. We omit some technical details and
   assume that $r=2$.  

 Step 1. Iterations of the lemma on distribution of doubling index are used to show that
 if a harmonic function $u$ has a doubling index in cube $Q$ smaller than some large number $N$ and 
 we partition $Q$ into $K^3$ equal subcubes with $K \leq N$, then the number of  $Q_j$ such that
\begin{equation}
 N(Q_j) \geq N / 2^{c_1\log K /\log\log K} 
\end{equation}
is smaller than $K^{2-c_2}$, where $c_1,c_2$ are small positive numerical constants.

 Step 2. Assume that $B$ is a ball of radius 2, $\max_{\frac{1}{4}B} |u|=1$ and $|u(0,0,1)|=\max_{\frac{1}{2}B} |u|$.
 The stable growth assumption implies  $$  \max_{\frac{1}{2}B} |u| \geq 2^N  \textup{ and } \max_{2B} |u| \leq 2^{CN}. $$
  We consider cubes forming a lattice with sides parallel to $x,y,z$ axes and with side length $1/ \sqrt N$. 
 Let us denote by $Q_j$ those cubes from the lattice, which intersect $B$. The total number of $Q_j$ is comparable  to $(\sqrt N)^3$.

 To prove the key lemma it is enough to show that there are at least $N 2^{c\log N/ \log\log N}$ cubes $Q_j$ that contain a zero point.

It follows from the step 1 that most of $Q_j$ (all except probably $N^{1-c_5}$ cubes) satisfy  
 \begin{equation} \label{eq: good}
N_u(Q_j) \leq N / W,
 \end{equation}
  where $W= 2^{c_3\log N/\log\log N}$.
 We call $Q_j$ good if \eqref{eq: good} holds. 

\textbf{Claim 2.} If $Q_j$ and $Q_k$ are adjacent and good, then  $$\max_{\frac{1}{2}Q_j}|u|\leq \max_{\frac{1}{2}Q_k}|u| 2^{CN/W}.$$

 We split cubes $Q_j$ into groups so that the centers of the cubes in each group lie on a line parallel to $z$-axis. We call such groups of cubes  tunnels.
Each tunnel has at most $C\sqrt N$ cubes.
 
 \textbf{Observation}. By step 1 most of the tunnels have only good cubes. 

 Now, consider tunnels that contain only good cubes and at least one cube $Q_j$ with distance to the maximum point $(0,0,1)$  smaller than $1/\log^2 N$.
 It follows from Step 1 that the total number of such tunnels  is at least $c_5 N/\log^4 N$.

The proof of the key lemma is completed by the next proposition. 

\textbf{Proposition.} Assume that a tunnel $T$  contains only good cubes and at least one cube  with distance to the maximum point $(0,0,1)$ smaller than $1/\log^2 N$.
Then $T$ contains at least $c_6 W$ cubes $Q_j$ with zeroes of $u$.

Step 3 (Proof of proposition).  

 Since $T$ is parallel to $z$ axis it  contains at least one cube $Q_a$ in $\frac{1}{4}B$ and therefore
 there is at least one cube in $T$ with $\max_{Q_a}|u|\leq 1$. We also know that $T$ contains at least one cube $Q_b$ with distance to the maximum point $(0,0,1)$ smaller than $1/\log^2 N$ and it appears that $$\max_{\frac{1}{2}Q_b}|u| \geq 2^{c_4 N}.$$

 The proof of the latter statement is using the assumption of stable growth. 
We split the latter statement into several claims, but omit some of the details.

\textbf{Claim 3.} If $\rho \in (\frac{1}{\log^{100} N}, 1/8)$, then $$N_u(B_\rho(0,0,1)) \leq C \rho N.$$

 Claim 3 says that the doubling index with the center at the maximum becomes smaller when we decrease the radius. 
 The next claim gives a lower bound on the maximum in small balls near the maximum.

\textbf{Claim 4.} Assume that $\rho \in (\frac{1}{\log^{100} N}, 1/8)$ and  the distance from  a point $x$ to the maximum point $(0,0,1)$ is smaller than $\rho$.
  If a positive number $s \leq \rho/2$, then 
  $$ N_u(B_s(x))\leq C \rho N $$
and  
  $$ \max\limits_{B_s(x)} |u| \geq |u(0,0,1)|2^{-C\rho N  \log (\frac{\rho}{s})}.$$

\textbf{Claim 5.} If a distance from a cube $Q_b$  to the maximum point $(0,0,1)$ is smaller than $1/\log^2 N$, then 
$$\max_{\frac{1}{2}Q_b}|u|\geq |u(0,0,1)| 2^{-C N /\log N}$$
and therefore $\max_{\frac{1}{2}Q_b}|u| \geq 2^{c_4 N} $.

 Now we are ready to finish the proof of the proposition. We start going along $T$ from $Q_a$  with $\max_{\frac{1}{2}Q_a}|u|\leq 1$ to the cube $Q_b$ with
$\max_{\frac{1}{2}Q_b}|u| \geq 2^{c_4 N}$ and watch  how the maximum over cubes changes. The total multiplicative increment is at least $ 2^{c_4 N}$.  We consider twice smaller cubes because further we will apply the Harnack inequality.

If we have two adjacent cubes $Q_j$ and $Q_k$ where $u$ does not have any zeroes, then the Harnack inequality guarantees that 
$$\max_{\frac{1}{2}Q_j}|u|\leq C \max_{\frac{1}{2}Q_k}|u| .$$
The total number of cubes in $T$ is smaller than $C\sqrt N$ and therefore  the total multiplicative increment over pairs of adjacent cubes with no zeroes
 is smaller $C^{C\sqrt N}$, which is a way smaller than $ 2^{c_4 N}$. In particular the Harnack inequality guarantess that there is at least one cube in $T$
 with a zero of $u$. So the major part of the multiplicative increment comes from the pairs of cubes, where at least one of the cubes has a zero point. 

All cubes in $T$ are good and by Claim 2 for any adjacent cubes in $T$
 $$\max_{\frac{1}{2}Q_j}|u|\leq \max_{\frac{1}{2}Q_k}|u| 2^{CN/W}.$$

 In particular, we cannot realize the increment $ 2^{c_4 N}$ passing only one pair of adjacent cubes, even if there is a zero there.
Moreover Claim 2 guarentees that the total number of cubes in $T$ with zeroes of $u$ is at least $c_6 W$.

\section{In between real-analytic and smooth cases. }
Between February 13 and February 17, 2017 there was a workshop   \cite{IAS} at the Institute for Advanced Study on Emerging Topics: Nodal sets of
Eigenfunctions. The first author was giving a talk on the joint result with the second author and N.Nadirashvili:

\textbf{Theorem}. If $\Omega$ is a bounded domain in $\R^n$ with a smooth boundary and $\varphi_\lambda$ is any Laplace eigenfunction of $\Omega$ with Dirichlet boundary conditions, then
 $\h^{n-1}(Z_{\varphi_\lambda}) \leq C_{\Omega} \sqrt \lambda \log (\lambda +e)$.

 During the talk Fedor Nazarov removed a half of the proof, which appeared to be unnecessary, simplified the argument and improved the bound to the optimal one:
 $$\h^{n-1}(Z_{\varphi_\lambda}) \leq C_{\Omega} \sqrt \lambda.$$
 
  The opposite inequality
$$\h^{n-1}(Z_{\varphi_\lambda}) \geq c_{\Omega} \sqrt \lambda.$$  
is also true (if we include the boundary of $\Omega$ in the nodal set  $Z_{\varphi_\lambda}$) and now it has  two different proofs: one is due to Donnelly and Fefferman and another proof involves the solution of Nadirashvili's conjecture.

\bibliographystyle{plain}  
\bibliography{nodal}

\end{document}